\documentstyle[theorem]{article}

\theorembodyfont{\rm}
\newtheorem{sdef}{Definition}[section]
\newtheorem{sprop}[sdef]{Proposition}
\newtheorem{sthm}[sdef]{Theorem}
\newtheorem{slem}[sdef]{Lemma}
\newcommand{\Db}[2]{{{\bf \Delta}^#1_#2}}
\newcommand{\Sb}[2]{{{\bf \Sigma}^#1_#2}}
\newfont{\Euler}{msbm10}
\newfont{\Fraktur}{eufm10}
\newcommand{\CCC}{{\mbox{{\Euler C}}}}
\newcommand{\DDD}{{\mbox{{\Euler D}}}}
\newcommand{\LLL}{{\mbox{{\Euler L}}}}
\newcommand{\MMM}{{\mbox{{\Euler M}}}}
\newcommand{\PPP}{{\mbox{{\Euler P}}}}
\newcommand{\RRR}{{\mbox{{\Euler R}}}}
\newcommand{\SSS}{{\mbox{{\Euler S}}}}
\newcommand{\TTT}{{\mbox{{\Euler T}}}}
\newcommand{\VVV}{{\mbox{{\Euler V}}}}
\newcommand{\WWW}{{\mbox{{\Euler W}}}}
\newcommand{\frnT}{{\mbox{{\Fraktur T}}}}
\newcommand{\omlom}{\omega^{<\omega}}
\newcommand{\omom}{\omega^\omega}
\newcommand{\omuplom}{\omega^{\uparrow<\omega}}

\newcommand{\pred}{{\rm pred}}
\newcommand{\Succ}{{\rm Succ}}
\newcommand{\stem}{{\rm stem}}
\newcommand{\et}{\,{}\hat{}\,}
\newcommand{\la}{\langle}
\newcommand{\ra}{\rangle}
\newcommand{\ran}{{\rm ran}}
\newcommand{\dom}{{\rm dom}}
\newcommand{\amal}{{\rm amal}}
\newcommand{\qed}{\begin{flushright}{\sf q.e.d.}\end{flushright}}
\newcommand{\Ord}{{\rm Ord}}
\newcommand{\rk}{{\rm rk}}

\begin{document}
\begin{center}
{\LARGE Uniform Unfolding and Analytic Measurability}\\[0.5cm]
{\large Benedikt L\"owe}\footnote{The author expresses his deep
gratitude towards J\"org Brendle who introduced him into the matter and
offered countless pieces of advice.
The paper was completed whilst the author held a grant
by the Studienstiftung des Deutschen Volkes.}\\
Eberhard--Karls--Universit\"at T\"ubingen\\
{\footnotesize loewe@rhein.philosophie.uni-tuebingen.de}\\
Humboldt--Universit\"at zu Berlin\\
{\footnotesize loewe@wega.mathematik.hu-berlin.de}\\
\end{center}
\section{Introduction}

There are many natural forcing notions adding real numbers, most
notably the precursors of this class of forcing notions, namely
{\sc Cohen} and Random forcing.
Whereas Random forcing has been
deliberately constructed in connection with the $\sigma$--algebra of
{\sc Lebesgue} measurable sets and {\sc Cohen} forcing is represented
by the algebra of {\sc Borel} sets modulo the $\sigma$--ideal of meager
sets, the other classical forcing notions, which add real numbers,
e.g. {\sc Sacks} forcing, {\sc Miller} forcing and {\sc Laver} forcing
are not connected with a $\sigma$--algebra by virtue of their construction.
But as one can see easily, each of these classical forcing notions 
corresponds naturally to a $\sigma$--algebra\footnote{The early investigation
of these algebras and ideals began with \cite{Mar}.}. {\sc Miller} alludes
to this fact in the rather vague formulated problem (11.10)
of his problem list\footnote{\cite{M93}}. The main point of the problem is:
What kinds of well--known results can be proved for these
$\sigma$--algebras?\\
E.g. it is known that the analytic sets are {\sc Lebesgue} measurable and have
the {\sc Baire} property, and we have the first sets which do not have this
property on the $\Db{1}{2}$ level in {\bf L}.
So, given a forcing notion $\PPP$ adding real numbers, the following
two natural questions are connected with {\sc Miller}'s problem (11.10):
\begin{itemize}
\item Are the analytic sets $\PPP$--measurable, i.e. are they elements
of the $\sigma$--algebra naturally connected with $\PPP$?
\item Do we have $\Db{1}{2}$ counterexamples in {\bf L}? If yes, how can
we characterize the axiom candidate ``All $\Db{1}{2}$ sets are 
$\PPP$--measurable" (similarly for $\Sb{1}{2}$)?
\end{itemize}
In this paper we will tackle the first of these questions. The second
question is answered by {\sc Brendle} and the present author in
\cite{BL99} for {\sc Hechler} and {\sc Miller} forcing. We will not only
show that for all prominent forcing notions the analytic 
$\PPP$--measurability is provable in ZFC, but provide a uniform
proof technique for questions of analytic measurability based on the
so called "{\sc Solovay}'s Unfolding Trick". 

\section{Preliminaries}
\subsection{Notation and Definitions}
A {\it tree} is a subset of $\omlom$ closed under initial segments.\\
Let $T$ be a tree.
We call a real $r\in\omom$ a {\it branch} of $T$ if
$\forall n: r|_n\in T$. The set of all branches in $T$ is denoted by
$[T]$. The (immediate) predecessor of a node $t\in T$ is uniquely determined
and is denoted by $\pred(t)$.
Obviously, if 
$\sup(s\in\omlom:\forall x\in[T](s\subseteq x))$ does not exist, there is
a real $y_T\in\omom$ with 
$[T]=\{y_T\}$. 
We then call the following object the {\it stem} of $T$:
$$\stem(T) := \left\{
\begin{array}{rl}\sup(s\in\omlom:\forall x\in[T](s\subseteq x))
&\mbox{if the supremum exists}\\
y_T & \mbox{else}\end{array}\right.$$
To denote that a tree is a subtree with strictly longer stem
we write
$$T\ll T':\iff T\leq T'\mbox{ and }\stem(T)\neq\stem(T')$$
For a node 
$t\in T$ we denote by $\Succ(t) :=
\{s\in T : \exists n\in\omega(t\et\la n\ra=s)\}$ the set of its
(immediate) successors.\\
For a tree $T$ and a finite sequence $s\in\omlom$ we define
$T\uparrow s:=\{t\in T:s\subseteq t\mbox{ or } t\subseteq s\}$.
A node of a tree $T$ is called 
{\it splitting node}, if it has more than one successor, and
$\omega$--{\it splitting node}, if it has infinitely many successors.\\[0.5cm]
We call a notion of forcing $\PPP$ {\it arboreal}\footnote{This notion
is far more general than the ``perfect tree property''
of \cite{GJ91}, but in the applications we have to construct
subtrees with specific properties and we need for these constructions
many of the properties of ``perfect tree forcing''. In particular
all of the investigated forcing notions have the perfect tree property.}
if there is a set $\frnT$
of trees partially ordered by inclusion, $\PPP$ is order--isomorphic
to $\frnT$ and $\frnT$ has the following property:
$$\forall T\in\frnT\forall N>|\stem(T)|\exists T'\leq T:|\stem(T')|\geq N$$
If $\PPP$ is an arboreal forcing, we identify $p\in\PPP$ with the
tree $T_p$.\\
Without loss of generality all arboreal forcings have a largest element
{\bf 1} with 
$[{\bf 1}]=\omom$.\\

\subsection{{\it dramatis personae}}
\begin{enumerate}
\item A tree $L\subseteq\omlom$ is called {\sc Laver} tree, if all
nodes above the stem are $\omega$--splitting nodes.
We call the set of all {\sc Laver} trees ordered by inclusion
{\sc Laver} forcing $\LLL$.
\item A tree 
$M\subseteq\omlom$ is called {\it superperfect}, if
every splitting node is an $\omega$--splitting node and every node
has a successor which is a splitting node (and therefore an 
$\omega$--splitting node).
{\sc Miller} forcing $\MMM$ is the set of all superperfect trees
ordered by inclusion.
\item In analogy to the definition of $\MMM$ we call a tree
$P\subseteq 2^{<\omega}$ {\it perfect}, 
if below every node there is a splitting node and define {\sc Sacks}
forcing $\SSS$ to be the set of all perfect trees ordered by inclusion.
\item We call a perfect tree $P$ {\it uniform}, 
if it has the following property:\\
If $t_1,t_2\in P$ with $|t_1|=|t_2|$, then $t_1\et\la 0\ra\in P
\iff t_2\et\la 0\ra\in P$ and $t_1\et\la 1\ra\in P
\iff t_2\et\la 1\ra\in P$.\\
The set of all uniform perfect trees ordered by inclusion is called
{\sc P\v r\'ikr\'y--Silver} forcing $\VVV$.\\
The uniformity of $\VVV$ obliges us to introduce the combinatorial
technique of amalgamation:
Let $T$ be a uniform perfect tree and $t_1,t_2\in T$ with $|t_1|=|t_2|$.
Define $R_1:=T\uparrow t_1$ and $R_2:=T\uparrow t_2$. If we have
$Q\leq R_2$, then we set
\begin{eqnarray*}
\amal(R_1, Q)&:=&\{t\in\omlom:t(i)=t_1(i)\mbox{ for }i<|t_1|\\
&&\mbox{ and }
\exists r\in Q:t(i)=r(i)\mbox{ for }i\geq|t_1|\}
\end{eqnarray*}
Obviously we construct a copy of $Q$ in $R_1$ so that 
$\amal(R_1,Q)\leq R_1$. 
\item Take the following set 
$\TTT:=\{\la s,A\ra\in\omuplom\times[\omuplom]^\omega:$ there is an
enumeration of $s\cup A$ with $a_0=s$ and ${\sf min}(\ran(a_{i+1}))
>{\sf max}(\ran(a_i))$ for all $i\in\omega\}$
and order it via
\begin{eqnarray*}
\la s,A\ra\leq\la t,B\ra&\Leftrightarrow& s\supseteq t, \forall a\in A\exists 
\{b_1,\dots, b_n\}\subseteq B (a=b_1\et\dots\et b_n),\\
&&\exists \{\beta_1,\dots,\beta_m\}\subseteq B
(\ran(s)\setminus\ran(t)=\ran(\beta_1\et\dots\et\beta_m))
\end{eqnarray*}
We call this partial ordering {\sc Matet} forcing\footnote{{\it cf.} 
\cite{Matet}}.
We define
$$x\in[\la s,A\ra]\iff s\subset x\wedge \exists A_0=\{c_1,c_2,\dots\}
\in [A]^\omega
\mbox{, so that } x=s\et c_1\et c_2\et\dots$$
and let $T_{\la s,A\ra}$ be the tree of all finite initial sequences
of reals in $[s,A]:=[\la s,A\ra]$.\\
As for {\sc P{\v r}\'ikr\'y--Silver} forcing we define for
$T:=\la s,A\ra$ and $S:=\la t, B\ra$ with $\la s, B\ra\leq \la s, A\ra$:
$$\amal(T, S) := \la s, B\ra$$
\item We define Willowtree forcing $\WWW$\footnote{{\it cf.}
\cite{B95} for the reasons of introducing this forcing and other
information about the connections between the mentioned forcing
notions.} by
\begin{eqnarray*}
\la f,A\ra\in\WWW&\iff& A\in[[\omega]^{<\omega}]^\omega,\forall a,a'\in A:
\end{eqnarray*}
We call the elements of $\WWW$ willowtrees.
$\WWW$ is arboreal via
\begin{eqnarray*}
T_W&:=&\{s\in\omlom: \forall n<|s|(n\in\dom(f)\to s(n)=f(n))\\
&&\hskip0.5cm\wedge\,\forall a\in A(|\ran(s|_a)|=1)\}
\end{eqnarray*}
where $W=\la f,A\ra$.\\
For every willowtree
$W=\la f,A\ra$ we order the elements of $A$ as follows:
\[i \leq j ~\iff ~\min(a_i)\leq \min(a_j)\]
Since the elements of $A$ are disjoint sets, this is a linear ordering.
Define the following relation on $\WWW$:
\[W\geq_i W'~:\iff~\forall j\leq i:a_j\in A'\]
\begin{sprop}
When $W^0\geq_0 W^1\geq_1 W^2\geq_2\dots$ is a sequence of willowtrees,
then $\la W^i:i\in\omega\ra$ is a fusion sequence.
\end{sprop}
Because of the affinity of Willowtree forcing with the uniform
forcing notions $\VVV$ and $\TTT$ we need a form of amalgamation:
For every willowtree $W$,
$\sigma\in 2^{<\omega}$ and $i=|s|$, set
\[W^\sigma_i := \la f\cup\bigcup_{j=0}^{i-1} (a_j\times \{\sigma(j)\}),
A\setminus\{a_0,\dots,a_{i-1}\}\ra=:\la f_i^\sigma, A^\sigma_i\ra\]
For amalgamation let now $W$ be an arbitrary willowtree,
$|\sigma|=|\sigma'|=i$,$V\leq W_i^\sigma$ and $V = \la g,B\ra$.\\
In this case we define the function
\[h:=\left\{\begin{array}{rl}f^{\sigma'}_i & \mbox{ on }\bigcup_{j=0}^{i-1}
a_j\\
g&\mbox{ else}\end{array}\right.\]
and with this $\amal(W^{\sigma'}_i,V):=\la h,B\ra$.
\end{enumerate}

\subsection{Measurability}
We can connect every arboreal forcing naturally to a notion of measurability.
In addition to the forcings defined there are the well-known
forcing notions of {\sc Cohen} forcing $\CCC$, {\sc Hechler} forcing
$\DDD$ and {\sc Mathias} forcing $\RRR$\footnote{For definitions {\it cf. e.g.}
\cite{J84}, \cite{Truss}}.
These forcings form topology bases for the {\sc Baire} topology ${\cal C}$, the
dominating topology ${\cal D}$ and the {\sc Ellentuck} topology 
${\cal R}$ respectively. The forcings are therefore
quite naturally connected to the $\sigma$--algebra of sets with the
{\sc Baire} property in these topologies. The ${\cal C}$--, ${\cal D}$--
and ${\cal R}$--meager sets are also called $\CCC$--, $\DDD$--
and $\RRR$--null sets.\\
In analogy to this situation we define in the case of non--topological 
forcings $\PPP$ a set of real numbers $A$ ($A\subseteq\omom$ or 
$A\subseteq 2^\omega$ according to the 
definition of $\PPP$) to be $\PPP$--measurable if 
\[\forall p\in\PPP\exists p'\leq p:([p']\cap A=\emptyset\mbox{ or }
[p']\cap \omom\setminus A=\emptyset)\]
\subsection{The {\sc Banach--Mazur} games}
The Polish school of set theorists and topologists commenced
the study of analytic sets via infinite games.
For this reason they defined the so--called {\sc Banach--Mazur} games
on a topological space.\footnote{For proofs of all results in this
section we refer 
to \cite{K95}, p.51 {\it sqq} and 149 {\it sqq}.}
\begin{sdef}
Let $X$ be a topological space and $A\subseteq X$. 
Then we state the rules of the 
{\sc Banach--Mazur} {\it game} $G_X(A)$:
\begin{itemize}
\item There are two players I and II, I begins.
\item The players play open sets $U^I_i$ and $U^{II}_i$ respectively in turn.
\item There are the following restrictions: $U^I_{i+1}\subseteq U^{II}_i$ and
$U^{II}_i\subseteq U^I_i$
\item After playing countably many turns, I wins if
$\emptyset\not =\bigcap_{i<\omega} U^I_i\subseteq A$. Else II wins.
\end{itemize}
\end{sdef}
\begin{sdef}
Let $X$ be a topological space and $A\subseteq X\times\omom$.
Then we state the rules of the 
{\sc Banach--Mazur} {\it game in the plane} $G^2_X(A)$:
\begin{itemize}
\item There are two players I and II, I begins.
\item The players play in turn; II plays open sets
$U^{II}_i$, I plays pairs of open sets and natural numbers
$\la U^I_i, n_i\ra$.
\item There are the following restrictions: $U^I_{i+1}\subseteq U^{II}_i$ 
and $U^{II}_i\subseteq U^I_i$
\item After playing countably many turns, I wins, if
$\emptyset\not=(\bigcap_{i<\omega} U^I_i)\times \{(n_i)_{i<\omega}\}
\subseteq A$. Else II wins.
\end{itemize}
\end{sdef}
We distinguish between the games $G_X(A)$ or $G^2_X(A)$ on the one hand,
which is the set of all legal sequences according to the rules mentioned
above, and a {\it run of the game} on the other,
which is one particular sequence following the rules.
A function on the initial segments of runs to the appropiate mathematical
objects to play in the next turn is called a {\it strategy}. Obviously we
call a strategy {\it winning} if it guarantees that its user will be
the winner of the game.
A set $A$ is called {\it determined} for a game, if either I or II has a 
winning strategy in this game and a topological space is called
{\sc Choquet}, if I has a winning strategy for $G_X(X)$.\\
We now state the fundamental theorems for {\sc Banach--Mazur} games
which decide our questions for the topological forcings $\CCC$, $\DDD$
and $\RRR$:
\begin{sthm}[{\sc Gale--Stewart} 1953]\label{GSOriginal}
Let $A$ be a closed set in $\omom$ or $(\omom)^2$.
Then $A$ is determined for the games $G_X(A)$ and $G^2_X(A)$.
\end{sthm}
\begin{sthm}[{\sc Banach--Mazur}]\label{BanachMazur}
Let $X$ be a {\sc Choquet} topological space which is a
refinement of a metric space and  $A\subseteq X$. Then for $G_X(A)$:
\begin{enumerate}
\item I has a winning strategy $\Rightarrow$ $A$ is comeager in a non--empty
open set
\item II has a winning strategy $\Rightarrow$ $A$ is meager
\end{enumerate}
\end{sthm}
\begin{sthm}[{\sc Solovay}s Unfolding Trick]\label{Solo}
Let $X$ be a {\sc Choquet} topological space, which is a refinement of a
metric space, $C\subseteq X\times\omom$ and $A$ the projection of $C$
on $X$. Then for $G^2_X(C)$:
\begin{enumerate}
\item I has a winning strategy $\Rightarrow$ $A$ is comeager in a non--empty
open set
\item II has a winning strategy $\Rightarrow$ $A$ is meager
\end{enumerate}
\end{sthm}
As a consequence of all these results we get in the topological case:
\begin{sthm}\label{SiTop}
All analytic sets are $\PPP$--measurable for $\PPP\in\{\CCC,\DDD,
\RRR\}$.\footnote{{\it cf.} \cite{M80}, p. 299{\it sq}}
\end{sthm}
\section{Generalized {\sc Banach--Mazur} games}
\subsection{Basic Notions and Generalized {\sc Gale--Stewart} theorem}
Our goal is now to prove a result similar to the
{\sc Banach--Mazur} theorem for the non--topological forcings.
We define the {\it generalized {\sc Banach--Mazur} game} and its
variant in the plane\footnote{There are some precursors of these
notions in \cite{KechrisForcing}, but they focus on topological
properties of forcing.}:
\begin{sdef}
Let $\PPP$ be an arboreal forcing and $A\subseteq\omom$.
Then the generalized {\sc Banach--Mazur} game $G_\PPP (A)$ is
played with the following rules:\\
I begins and I and II play in turn 
forcing conditions $p_i^I$ und $p_i^{II}$ with the following restricting
property:
$$p_{i+1}^{I}\ll p_{i}^{II}\ll p_i^{I}$$
Then $f:=\bigcap [p_i^I]$ is obviously a real number.
I wins, if $f\in A$, otherwise II wins.
\end{sdef}
\begin{sdef}
Let $\PPP$ be an arboreal forcing and $A\in(\omom)^2$.
Then the generalized {\sc Banach--Mazur} game in the plane
$G^2_\PPP (A)$ is played with the following rules:\\
I begins and plays pairs of forcing conditions $p_i^I$ 
and natural numbers 
$n_i$, 
and II plays forcing conditions
$p_i^{II}$ with the following property:
$$p_{i+1}^{I}\ll p_{i}^{II}\ll p_i^{I}$$
Then $f:=\bigcap [p_i^I]$ is obviously a real number.
I wins, if $\la f,\la n_i:i\in\omega\ra\ra\in A$
and II wins, if $\la f,\la n_i:i\in\omega\ra\ra\not\in A$.
\end{sdef}
To get a result analogous to the topological case we have to
generalize the {\sc Gale--Stewart} theorem:
For a set $A$ and a forcing condition $p$ define the relative game
$G_\PPP(A,p)$ or $G^2_\PPP(A,p)$ simply by postulating that
all played conditions lie below $p$.
\begin{sprop}\label{GSAllg}
If $A\subseteq\omom$ is a closed set, $\PPP$ an arboreal forcing
and $p\in\PPP$, then all games $G_\PPP(A,p)$, $G_\PPP(A)$,
$G^2_\PPP(A,p)$ and $G^2_\PPP(A)$ are determined.
\end{sprop}
{\bf Proof :}\\
One can easily see: if II has no winning strategy for
$G_\PPP(A,p)$ (or $G^2_\PPP(A,p)$) then there is a $q\leq p$ with
the property that for all $r\leq q$ II has no winning strategy in
$G_\PPP(A,r)$.\\
Let $A$ be closed. Suppose that II has no winning strategy for 
$G_\PPP(A,p)$ so there is such a $q\leq p$. Then I can choose a $q$
with strictly longer stem. Regardless of what $q'$ II answers, II will
have no winning strategy in $G_\PPP(A, q')$, so we can iteratively define
a strategy for I.\\
We still have to show that this strategy is winning or, in other
words, that $\bigcup_{n\in\omega} \stem(p_n)=\bigcap_{n\in\omega} 
[p_n]=:f\in A$.\\
Suppose not, then one can find a finite sequence $s\subseteq f$
such that $[s]\subseteq \omega^\omega\setminus A$, because the complement
of $A$ is open. Without loss of generality we have for some $n_0$:
$s=\stem (p_{n_0})$.
But then II would have a trivial winning strategy for $G_\PPP(A,p_{n_0})$ 
in contradiction to our assumption.\\
For the games in the plane we regard in the proof simply the product
ordering in $\PPP\times\omlom$ instead of $\PPP$.\qed
For the general context of uniform unfolding we need two fundamental
notions:\\
\begin{sdef}\label{K}
Let $\PPP$ be an arboreal forcing and $\la\tau_i:i\in\omega\ra$ a sequence of
strategies for one fixed player in the game $G_\PPP$.
A partial function
$\alpha:\PPP\to\PPP$
will be called a $\PPP$--{\it strategic fusion} for $\la\tau_i:i\in
\omega\ra$ if $\alpha(T)$ has the following properties:
\begin{itemize}
\item[(K1)] $\alpha(T)\leq T$
\item[(K2)] There is a function $S:[\alpha(T)]\times \omega\to
\{P\in\PPP:P\leq T\}:\la x,b\ra\mapsto S_{x,b}$ with 
$$S_{x,b+1}\ll \tau_b(S_{x,b})\ll S_{x,b}$$
\item[(K3)] $\dom(\alpha)=\PPP$ if $\tau_i$ is a sequence of strategies
for player II and $\dom(\alpha)=\{\tau_0({\bf 1})\}$,
if $\tau_i$ is a sequence of strategies for player I.
\end{itemize}
\end{sdef}
\vskip1cm
In most cases we can even construct a sequence $\la T_\sigma : \sigma\in
\omlom\ra$ with $\alpha(T)=\bigcap_{i\in\omega}\bigcup_{|\sigma|=i}T_\sigma$
(hence the name fusion) with the properties:
\begin{itemize}
\item[(C1)] For each $x\in[\alpha(T)]$ there is an increasing sequence
$\la\sigma_i : i\in\omega\ra$ in $\omlom$ with $|\sigma_i|=i$,
so that $x\in[T_{\sigma_i}]$ for all $i\in\omega$.
\item[(C2)] For each $\sigma\in\omlom$ there is an $S\in\PPP$ such that
$$T_{\sigma}\ll \tau_{|\pred(\sigma)|}(S)\ll S\leq T_{\pred(\sigma)}$$
\end{itemize}
In this case we call the function $f$ {\it constructive}.\\
To see that such an $f$ really is a strategic fusion, let $x\in [f(T)]$.
Because of (C1) there is a sequence $\la\sigma_i:i\in\omega\ra$ so that
$x\in[T_{\sigma_i}]$ for all $i\in\omega$. 
Take the tree $S$ whose existence is postulated in (C2) to be
$S_{x,b}$. Obviously the so defined function $S$ has property (K2).
\begin{sdef}
We say that an arboreal forcing $\PPP$ has the {\it linear
dichotomy property} if for every $A\subseteq\omom$:
\begin{enumerate}
\item If I has a winning strategy in $G_\PPP(A)$ then there is a $q\in\PPP$,
so that $[q]\subseteq A$
\item If II has a winning strategy in $G_\PPP(A)$, then there is for every
$p\in\PPP$ a $q\leq p$ with $[q]\cap A=\emptyset$
\end{enumerate}
\end{sdef}
\begin{sdef}
We say an arboreal forcing has the {\it planar dichotomy property}
if for every set
$C\subseteq(\omom)^2$ and its projection $A$ the following hold:
\begin{enumerate}
\item If I has a winning strategy in $G^2_\PPP(C)$, then there is a
$q\in\PPP$, so that $[q]\subseteq A$
\item If II has a winning strategy in $G^2_\PPP(C)$, then there is for every
$p\in\PPP$ a $q\leq p$ with $[q]\cap A=\emptyset$
\end{enumerate}
\end{sdef}
\subsection{The Unfolding Theorem}
\begin{sthm}\label{Eindim}
If there is for every strategy $\tau$ of an arbitrary player a
$\PPP$--strategic fusion $\alpha$ for the constant sequence 
$\la\tau:i\in\omega\ra$,
then $\PPP$ has the linear dichotomy property.
\end{sthm}
{\bf Proof :}\\
Let $\tau$ be a winning strategy for II and $P\in\PPP$
arbitrary. Let $Q:=\alpha(P)\leq P$ according to (K1).
We have to show that $[Q]\cap A=\emptyset$. Take $f\in[Q]$ 
arbitrary. Then we have with (K2) for each $i\in\omega$:
$$S_{f,i+1}\leq \tau_i(S_{f,i})=\tau(S_{f,i})\leq S_{f,i}$$
Obviously the sequence
$$S_{f,0}, \tau(S_{f,0}), S_{f,1}, \tau(S_{f,1}),
\dots$$
is a run of the game $G_\PPP(A)$ according to $\tau$. Therefore $f\notin A$.\\
On the other hand if $\tau$ is a winning strategy for I, so we can construct
the same sequence with one difference: (K3) says that $\alpha$ is
not defined for arbitrary $p\in\PPP$, so the constructed run has to
begin with the initial value of $\tau$, that is $\tau({\bf 1})$.
This accounts for the asymmetry between I and II in the definition
of the dichotomy property.\qed
\begin{sthm}[Unfolding Theorem]\label{ZweidimDich}
If there is a $\PPP$--strategic fusion for any sequence of strategies for any
of the two players then $\PPP$ has the planar dichotomy property.
\end{sthm}
{\bf Proof :}\\
Let $\beta$ be a bijection between $\omlom$ and $\omega$ having the property:
$s\subseteq t\Rightarrow \beta(s)\leq\beta(t)$.\\
Firstly, regard a winning strategy $\tau$ for II.
We define for $P\in\PPP$ the following sequence of strategies:
$\tau_i(P):=\tau(\la P,\beta^{-1}(i)\ra)$.
If $T\in\PPP$ is arbitrary, then there is according to our assumption
the strategic fusion $\alpha$, so
$T':= \alpha(T)\leq T$ via (K1).\\
We claim that $[T']\cap A=\emptyset$, i.e., more precisely:
If $x\in [T']$, then for all $y\in\omom$ it holds that $\la x,y\ra\notin C$.\\
Let now be $x\in [T']$ and $y\in\omom$. Then the sequence
$(y|_k)_{k\in\omega}\subseteq \omlom$ is increasing in
$\omlom$, and therefore $b_k := \beta(y|_k)$ is (because of the property
postulated for $\beta$) an infinite increasing sequence in $\omega$.
Via (K2) we get 
$$S_{x,b_{k+1}}\leq \tau_{b_k}(S_{x,b_k})\leq S_{x,b_k}$$
and so 
$$S_{x,b_0}\gg \tau_{b_0}(S_{x,b_0})\gg S_{x,b_1}\gg \tau_{b_1}
(S_{x,b_1})$$ is a decreasing sequence of trees.
So we can define a run
\begin{eqnarray*}
\la S_{x,b_0},y|_0\ra\geq\la\tau(\la S_{x,b_0},y|_0\ra),y|_0\ra&\geq&
\la S_{x,b_1},y|_1\ra\geq\la\tau(\la S_{x,b_1},y|_1\ra),y|_1\ra
\geq\dots\\
&\geq&\la S_{x,b_{k}},y|_{k}\ra\geq\la\tau(\la S_{x,b_{k}},y|_{k}
\ra),y|_{k}\ra\geq\dots
\end{eqnarray*}
of the planar game $G^2_\PPP(C)$ following the winning strategy
$\tau$. Therefore $\la x,y\ra\not\in C$ and because $y$ was arbitrary,
$x\not\in A$.\\
Secondly, let $\tau$ be a winning strategy for I. Now we have no
function from $\PPP\times\omlom$ to $\PPP$, but a function from 
$\PPP\times\omlom$ to $\PPP\times\omega$. So
$\tau(P,\beta^{-1}(i))$ has two components. We denote the first of them
by $\tau_i(P)$, and the second by $n_i(P)$.\\
By the assumption, we have a strategic fusion $\alpha$ for the sequence
$\la\tau_i:i\in\omega\ra$ so set $T':=\alpha(\tau_0({\bf 1}))$.
We now have to show that $[T']\subseteq A$, i.e.
$\forall x\in[T']\exists y\in\omom:\la x,y\ra\in C$.
Take now $x\in [T']$ and construct by recursion:
$$T_0 := \tau_0 ({\bf 1})$$
$$\sigma_0 := \la n_0 ({\bf 1})\ra$$
$$b_1 := \beta(\sigma_0)$$
$$T_1 := S_{x,b_1}\leq T_0$$
Then we have $S_{x,b_1+1}\leq\tau_{b_1}(T_1)\leq T_1$.
$$\sigma_{i+1}:= \sigma_i\et \la n_{b_{i+1}}(T_{i+1})\ra$$
$$b_{i+2}:=\beta(\sigma_{i+1})$$
$$T_{i+2}:= S_{x,b_{i+2}}$$
By the assumption $b_{i+1}\geq b_i+1$ and so $T_{i+1}\leq S_{x,b_i+1}\leq
\tau_{b_{i}}(T_i)\leq T_i$. So we
get the following run of the game $G^2_\PPP$:
\begin{eqnarray*}
\la T_0, \sigma_0\ra =  \la \tau_0({\bf 1}),\la n_0({\bf 1})\ra\ra  =
\tau(\la{\bf 1},\la\ra\ra)
&\geq &\la T_1, \sigma_0\ra\geq \la \tau_{b_1}(T_1), \sigma_1\ra\\
&\geq & \dots\\
&\geq & \la T_i, \sigma_{i-1}\ra \geq \la \tau_{b_i}(T_i), 
\sigma_i\ra\\
&\geq & \dots
\end{eqnarray*}
which is a run according to $\tau$ having $x$ in the first component as
a result and we have for $y:=\bigcup_{i\in\omega} \sigma_i$:
$$\la x,y\ra\in C$$\qed

\subsection{Strategic Fusions}\label{Spielverlauf}
After having reduced the proof of uniform unfolding to the
existence of strategic fusion by \ref{ZweidimDich}, we have to
give explicitly a strategic fusion for all the forcings
mentioned above.
\paragraph{{\sc Sacks} forcing}
For a perfect tree $P$ we denote by $h_P$ the first splitting node.
Let $T$ be a fixed perfect tree and $\la\tau_i:i\in\omega\ra$ a 
sequence of strategies.
Regard the following sequence of perfect trees.
\[P_{\la\ra} := T\]
\[P_{\sigma\et\la 0\ra} := \tilde P_\sigma\uparrow h_{\tilde P_\sigma}\et\la 
0\ra\]
\[P_{\sigma\et\la 1\ra} := \tilde P_\sigma\uparrow h_{\tilde P_\sigma}\et\la 
1\ra\]
\[\tilde P_\sigma := \tau_{|\sigma|} ( P_\sigma)\]
Because the $P^{(i)}:=\bigcup_{|\sigma|=i} \tilde 
P_\sigma$ are a fusion sequence we can define $\alpha(P):=
\bigcap_{i\in\omega}P^{(i)}$ and $\alpha$ is even constructive.

\paragraph{{\sc Miller} forcing}
The strategic fusion is exactly analogous to the {\sc Sacks} run
construction. Simply substitute 2 by $\omega$.

\paragraph{{\sc Laver} forcing}
If $\tau$ is a strategy for the {\sc Laver} game and $T$ is a 
{\sc Laver} tree, so we define a rank function:
\[\rk^\tau_T : T\to \Ord\cup\{\infty\}\]
\[\rk^\tau_T (t) = 0~:\iff~\exists S\ll T\,:\,\stem(\tau(S))=t\]
\[\rk^\tau_T(t)\leq\alpha~:\iff~\forall^\infty t'\in \Succ(t)\,:\,
\rk^\tau_T(t')<\alpha\]
\[\rk^\tau_T(t) =\infty~~~\mbox{othwerwise}\]
Then we can easily prove the following
\begin{slem}
If $\tau$ is a strategy, we have for every node $t\in T$:
$\rk^\tau_T(t)\in\Ord$.
\end{slem}
Assisted by this lemma we now define a strategic fusion:
Let $T$ be an arbitrary {\sc Laver} tree and $\la\tau_i:i<\omega\ra$ a
sequence of strategies. At first, we define recursively a sequence
of {\sc Laver} trees $\tilde H_\sigma$.
To begin the recursion we define  $H_{\la\ra}:=T$, $\tilde H_{\la\ra}:=T$,
$R_{\la\ra}:=T$ and $i_{\la\ra}:=0$.\\
Let $\tilde H_\sigma$, $R_\sigma$ and $i_\sigma$ be already defined and
suppose that we have $\rk^{\tau_{i_\sigma}}_{R_\sigma}
(\stem(\tilde H_\sigma))\neq 0$ (we will show after the construction
that this is always the case). Then $\stem(\tilde H_\sigma)$ has 
infinitely many successors with smaller rank. We enumerate these 
successors by $\zeta_k$ and define recursively:
$$H_{\sigma\et\la n\ra} := \tilde H_\sigma\uparrow\zeta_n$$
$$\tilde H_{\sigma\et\la n\ra}:=\left\{\begin{array}{lcr}
\mbox{{\bf Case 1 : }}& H_{\sigma\et\la n\ra} & \mbox{if } 
\rk^{\tau_{i_\sigma}}_{R_\sigma}
(\stem(H_{\sigma\et\la n\ra}))>0\\
\mbox{{\bf Case 2 : }}& \tau_{i_\sigma}(S) & \mbox{if } 
\rk^{\tau_{i_\sigma}}_{R_\sigma}(\stem(H_{\sigma\et\la n\ra}))=0
\end{array}\right.$$
In this definition, $S$ is the condition below $R_\sigma$
with $\stem(H_{\sigma\et\la n\ra})=
\stem(\tau_{i_\sigma}(S))$ (which exists according to the definition of
the rank).\\
If we have Case 1 for $n\in\omega$, then we define:
$$R_{\sigma\et\la n\ra}:= R_\sigma$$
$$i_{\sigma\et\la n\ra}:= i_\sigma$$
If we have Case 2, we define instead of that:
$$R_{\sigma\et\la n\ra}:= \tilde H_{\sigma\et\la n\ra}$$
$$i_{\sigma\et\la n\ra}:= i_\sigma+1$$
If we have Case 2 at a tree, we call this tree a {\it switching point}.
If we have $k$ switching points among the predecessors of another
switching point, we call it a {\it switching point of order $k+1$}.
As we remarked above, 
$\rk^{\tau_{i_\sigma}}_{R_\sigma}(\stem(\tilde H_{\sigma}))$ can never be 
zero, because either $\tilde H_\sigma$ was a switching point,
then we have $\tilde H_\sigma = R_{\sigma}$ and (since strategies always
strictly prolongate the stem) the image of a strategy cannot have the same
stem, or $\tilde H_\sigma$ was no switching point and the rank is
larger than zero.\\
As the ranks are descending chains of ordinals, we know that after each tree
there is a switching point. So the set of all switching points is
order isomorphic to $\omlom$ and also the set of the corresponding trees
$S$. Now take this set of all these trees $S$
and denote them by $\la T_\sigma : \sigma
\in\omlom\ra$. 
Then the map $\alpha(T) := \bigcap_{i\in\omega} T^{(i)}$ with
$T^{(i)}:=\bigcup_{|\sigma|=i} T_\sigma$ 
is a constructive strategic fusion.

\paragraph{{\sc P\v r\'ikr\'y--Silver} forcing}
As in the case of {\sc Sacks} forcing we denote the first splitting node
of $P$ by $h_P$.
We identify finite sequences from the set $2^n$ with the corresponding
binary numbers. When we refer to the binary numbers, we write an
upper index $^{[n]}$ to indicate the length of the original sequence.
So we have $P_{\la 1,0,0\ra}=P^{[3]}_1$ and
$P_{\la 1,0\ra}=P^{[2]}_1$.\\
Take the following sequence of uniform perfect trees for an arbitrary tree
$T$:
\[P_{\la\ra} := T\]
\[P_{\sigma\et\la 0\ra} := \tilde P_\sigma\uparrow h_{\tilde 
P_\sigma}\et\la 0\ra\]
\[P_{\sigma\et\la 1\ra} := \tilde P_\sigma\uparrow h_{\tilde 
P_\sigma}\et\la 1\ra\]
For $|\sigma|=n$ we have $2^n$ trees $P_\sigma$, which are ordered via the
natural ordering of the binary numbers corresponding to the $\sigma$s.
We define iteratively:
\[U^{[n]}_0:= P^{[n]}_0\]
\[T^{[n]}_i := \tau_n(U^{[n]}_i)\]
\[U^{[n]}_i := \amal(P^{[n]}_i,T^{[n]}_{i-1})\]
\[\tilde P^{[n]}_i := \amal(P^{[n]}_i, T^{[n]}_{2^n-1})\]
With the same arguments as for {\sc Sacks} forcing 
$\alpha(T):=\bigcap_{i\in\omega}\bigcup_{|\sigma|=i}
\tilde P_\sigma$ is a constructive strategic fusion.

\paragraph{Willowtree forcing}
Let $W:=\la f,A\ra$ be an arbitrary willowtree, where $A=\{a_i : i\in\omega\}$.
We define the following sequence of willowtrees:
\[T_{\la\ra} := W\]
Suppose that $T_\sigma$ is defined for $\sigma\in 2^{<\omega}$ with
$|\sigma|=n-1$ and write $T_\sigma=\la f_\sigma, A_\sigma\ra$. Then:
\[H_{\sigma\et\la 0\ra}:=\la f_\sigma\cup (a_{|\sigma|}\times\{0\}),
A_\sigma\setminus\{a_{|\sigma|}\}\ra\]
\[H_{\sigma\et\la 1\ra}:=\la f_\sigma\cup (a_{|\sigma|}\times\{1\}),
A_\sigma\setminus\{a_{|\sigma|}\}\ra\]
\[\tilde H^{[n]}_0:=\tau_n(H^{[n]}_0)\]
\[\tilde H^{[n]}_{i+1}:=\tau_n(\amal(H^{[n]}_{i+1},\tilde H^{[n]}_i))\]
\[T^{[n]}_{2^n-1} :=\tilde H^{[n]}_{2^n-1}\]
\[T^{[n]}_{i-1}:=\amal(\tilde H^{[n]}_{i-1},T^{[n]}_{2^n-1})\]
As one can see easily, $T^{(n)}$ is a fusion sequence and therefore 
$W':=\bigcap_{n\in\omega} T^{(n)}$ is a willowtree. So $\alpha(W):= W'$
has all properties of a constructive strategic fusion.

\paragraph{{\sc Matet} forcing}
Let $T:=\la s, A\ra$ be a {\sc Matet} condition where we have
$A:=\{a_i:1\leq i<\omega\}$.
Let $\la \tau_i : i\in\omega
\ra$ be a sequence of strategies. Define the mapping 
$\Theta_i := \tau_i\circ\dots\circ\tau_0$.
Then regard the following sequence of {\sc Matet} conditions:
$$T_0 := T$$
$$s_{\la 0\ra}:= s\et a_1$$
$$s_{\la 1\ra} := s$$
$$H_{\la 0\ra} := T\uparrow s_{\la 0\ra}$$
$$H_{\la 1\ra} := T\uparrow s_{\la 1\ra}$$
$$\tilde H_{\la 0\ra} := \Theta_0 (H_{\la 0\ra})$$
$$\tilde H_{\la 1\ra} := \amal (H_{\la 1\ra}, \tilde H_{\la 0\ra})$$
$$T_1 := \tilde H_{\la 0\ra}\cup\tilde H_{\la 1\ra}$$
Now write $T_1$ as $\la s, A^{(1)}\ra$ where $A^{(1)}=
\{ a^{(1)}_i : i\in\omega\}$. Now we suppose that in the $i$th
step the sequences $s_\sigma$ with $|\sigma|=i-1$ are already
constructed and that we have the tree $T_{i-1} = \la s, A^{(i-1)}\ra$.
Define
$$s_{\sigma\et\la 0\ra} := s_\sigma\et a^{(i)}_i$$
$$s_{\sigma\et\la 1\ra} := s_\sigma$$
$$H_\sigma := T\uparrow s_\sigma$$
Again we identify binary numbers and $0$--$1$ sequences and set
$$\tilde H^{[i]}_0 := \Theta_i (H^{[i]}_0)$$
$$\tilde H^{[i]}_k := \left\{\begin{array}{ll}
\amal(H^{[i]}_k, \tilde H^{[i]}_{k-1})&\mbox{if $k$ is odd}\\
\Theta_i(\amal(H^{[i]}_k, \tilde H^{[i]}_{k-1}))&\mbox{if $k$ is even}
\end{array}\right.$$
$$\bar H^{[i]}_{2^i-1}:= \tilde H^{[i]}_{2^i-1}$$
$$\bar H^{[i]}_k := \amal(\tilde H^{[i]}_k, \bar H^{[i]}_{k+1})$$
$$T_i := \bigcup_{|\sigma|=i} \bar H_\sigma$$
Now we claim that $\alpha(T):=\bigcap_{i\in\omega} T_i$ is a strategic fusion.
Let $x\in\bigcap_{i\in\omega} [T_i]$ and $b\in\omega$. Then there is an 
increasing sequence $\la k_i : i\in\omega\ra$ with
$$x = s\et a^{(k_0)}_{k_0} \et a^{(k_1)}_{k_1} \et\dots$$
Define $S_{x,0}:=\amal(\tilde H^{[k_0]}_{\ell-1}, H^{k_0}_\ell)$ where
$\ell<2^{k_0}$ is the appropriate integer such that $x\in[\tilde 
H^{[k_0]}_\ell]$. For $b>0$ take $$k(b) :=
{\sf min}\{k_i : b\leq k_i\mbox{ and }\forall c<b (S_{x,c}\mbox{ was
not constructed on level }k_i)\}$$
Again take the corresponding 
$\ell < 2^{k(b)}$ and then define 
$$S_{x,b}:=\Theta_{b-1}(\amal(\tilde H^{[k(b)]}_{\ell-1},H^{[k(b)]}_\ell))$$
This function verifies that $\alpha$ is a strategic fusion.

\subsection{Results}
Concluding we get:
\begin{sthm}\label{Analytic}
For forcing notions $\PPP\in \{\SSS,\MMM,\TTT,\VVV,\WWW,\LLL\}$
all analytic sets are $\PPP$--measurable.
\end{sthm}
{\bf Proof :}\\
Together with the strategic fusions from the preceding section and
\ref{ZweidimDich} this proof is similar to \ref{SiTop}:\\
Let $\PPP$ be one of the mentioned forcings and $A$ be analytic in
$\omom$. 
Then define $\PPP_A:=\{P\in\PPP:\exists R\leq P\mbox{ with }
[R]\subseteq A\}$. We have to show that for all
$Q\notin\PPP_A$ we have : $\exists P'\leq Q:[P']\cap A=\emptyset$.
Take $Q\notin\PPP_A$, so we have no $R\leq Q$ with
$[R]\subseteq A$. But $A\cap[Q]$ is analytic, too, so we have a
closed set $C$ in the plane with the projection $A\cap[Q]$.
With Proposition \ref{GSAllg} $C$ is determined, according to
Section \ref{Spielverlauf} all of the forcings have strategic fusions
and with \ref{ZweidimDich} we get:
Either $\exists P:[P]\cap\omom\setminus (A\cap[Q])=\emptyset$ or
$\forall P\exists P'\leq P:[P']\cap (A\cap[Q])=\emptyset$. Since
$Q\notin\PPP_A$ the first case is impossible.
Choose now $P:=Q$, then we have a $P'\leq Q$ with 
$[P']\cap (A\cap[Q])=[P']\cap A =\emptyset$.\qed
What we have won are previously unknown results for $\VVV$, $\WWW$ and
$\TTT$ (the results for $\SSS$, $\MMM$ and $\LLL$ were already known
implicitly because of results connected with the asymmetric 
games\footnote{cf. \cite{D64}, \cite{K77} and \cite{GRSS95}}) and a 
uniform method for proving the same result in case a new arboreal 
forcing should appear.


\begin{thebibliography}{999}
\bibitem[{\sc Brendle} 1995]{B95} J\"org {\bf Brendle}, Strolling through
paradise, Fundamenta Mathematicae 148 (1995), p. 1--25
\bibitem[{\sc Brendle--L\"owe} 199?]{BL99} J\"org {\bf Brendle}, Benedikt
{\bf L\"owe}, $\Db{1}{2}$-- and $\Sb{1}{2}$--Levels of 
{\sc Hechler}, {\sc Laver} and {\sc Miller} measurability,
{\it in preparation}
\bibitem[{\sc Davis} 1964]{D64} M.{\bf Davis}, Infinite games of
perfect information, Annals of Mathematical Studies 52 (1964), p. 85--101
\bibitem[{\sc Ellentuck} 1974]{E74} Erik {\bf Ellentuck},
A new Proof that analytic sets are {\sc Ramsey}, Journal of
Symbolic Logic 39 (1974), p. 163--165
\bibitem[{\sc Goldstern--Repick\'y--Shelah--Spinas} 1995]{GRSS95}
Martin {\bf Goldstern}, Miroslav {\bf Repick\'y}, Saharon {\bf Shelah},
Otmar {\bf Spinas}, On tree ideals, Proceedings of the American
Mathematical Society 123 (1995), p. 1575--1581
\bibitem[{\sc Groszek--Jech} 1991]{GJ91} M.{\bf Groszek},
Thomas {\bf Jech}, Generalized Iteration of Forcing,
Transactions of the American Mathematical Society 324 (1991),
p. 1--26
\bibitem[{\sc Jech} 1986]{J84} Thomas {\bf Jech}, Multiple Forcing,
Cambridge 1986 $[$ Cambridge Tracts in Mathematics 88 $]$
\bibitem[{\sc Kechris} 1977]{K77} Alexander S. {\bf Kechris}, On a notion of
smallness of subsets of the {\sc Baire} space, Transactions of the
American Mathematical Society 229 (1977), p. 191--207
\bibitem[{\sc Kechris} 1978]{KechrisForcing} Alexander S.{\bf Kechris},
Forcing in Analysis, in: Gert H.M\"uller, Dana S.Scott ({\it eds.}),
Higher Set Theory, Proceedings Oberwolfach 1977, Heidelberg 1978
$[$Lecture Notes in Mathematics 669$]$, p.277--302
\bibitem[{\sc Kechris} 1995]{K95} Alexander S.{\bf Kechris},
Classical Descriptive Set Theory, New York 1995 $[$ Graduate Texts
in Mathematics 156 $]$
\bibitem[{\sc Marczewski} 1935]{Mar} Edward {\bf Marczewski}
{\small{\bf (Szpilrajn)}}, Sur une classe de fonctions de M.
Sierpi\'nski et la classe correspondante d'ensembles, Fundamenta
Mathematica 24 (1935), p. 17--34
\bibitem[{\sc Matet} 1988]{Matet} Pierre {\bf Matet}, Some
filters of partitions, Journal of Symbolic Logic 53 (1988),
p. 540--553
\bibitem[{\sc Miller} 1993]{M93} Arnold W. {\bf Miller}, Arnie Miller's 
Problem list, in: Haim Judah ({\it ed.}), Set Theory of the Reals,
Israel Mathematical Conference Proceedings Vol 6 (1993), p. 645--654
\bibitem[{\sc Moschovakis} 1980]{M80} Yiannis N. {\bf Moschovakis},
Descriptive Set Theory, Amsterdam 1980 $[$Studies in Logic
and the Foundation of Mathematics 100$]$
\bibitem[{\sc Solovay} 1970]{Sol70} Robert M. {\bf Solovay}, A model of
set--theory in which every set
of reals is {\sc Lebesgue} measurable, Annals of Mathematics
92 (1970), p. 1--56
\bibitem[{\sc Truss} 1977]{Truss} John {\bf Truss}, Set having
calibre $\aleph_1$, in: R.Gandy, M.Hyland, {\it (eds.)}, Logic Colloquium
76, Amsterdam 1977, S. 595--612
\end{thebibliography}
\end{document}